\documentclass[12pt]{article}%
\usepackage{amssymb}
\usepackage{amsfonts}
\usepackage{amsmath}
\usepackage{graphicx}%
\begin{document}

\begin{center}
{\Large On the imbalance lattice of path-length sequences of binary trees}

\bigskip

S. Foldes* and S. Radeleczki**

*Tampere University of Technology,\ Finland

**University of Miskolc, Hungary

\bigskip

\bigskip\textbf{Abstract}
\end{center}

\textit{The existence of greatest lower bounds in the imbalance order of
path-length sequences of binary trees is seen to be a consequence of a joint
monotonicity property of the greater and lower expension operations. Path
length sequences that are join-irreducible in the imbalance lattice are
characterized.}

\bigskip

\bigskip

\textbf{1} \textbf{Introduction, terminology, notation}

\textbf{\bigskip}

Generally the framework and terminology of Stott Parker and Prasad Ram
$\left[  \text{SPPR}\right]  $ is followed in what follows.

\bigskip

For any sequence $x=$ $(x_{1},...,x_{n})$ of real numbers, $n\geq1$, we write
$\exp x$ for the sequence $(2^{-x_{1}},...,2^{-x_{n}})$, we write $\Sigma x$
for the sum $x_{1}+...+x_{n}$, and $last$ $x$ (resp $first$ $x$) for $x_{n}$
(resp. $x_{1}$). We also write $Sx$ for the sequence of partial sums
$(x_{1},x_{1}+x_{2},...,x_{1}+...+x_{n}).$ The \textit{suffix length }$suf$
$x$ of $x$ is the largest integer $k\leq n$ such that the last $k$ components
of $x$ are equal. A sequence $x$ of non-negative integers is a
\textit{path-length sequence }if $x_{1}\leq...\leq x_{n}$ and $\Sigma\exp=1.$
By Kraft's Theorem, this means that there is binary tree whose root-to-leaf
paths are in order of increasing lengths from left to right, the lengths being
$x_{1},...,x_{n}$ in that order $\left[  \text{K}\right]  $. Path-length
sequences have even suffix length. Between path-length sequences $l$ and $h$
with the same number of components, $l$ is said to be \textit{more balanced}
than $h$, in symbols $l\trianglelefteq h$, if $S\exp l\leq S\exp$ in the
componentwise order of vectors. This defines a partial order relation on the
set of path-length sequences with a given number $n$\ of components.

For a path-length sequence $l=$ $(l_{1},...,l_{n})$, the \textit{expansion in
position} $i$ is defined as thepath-length sequence $(l_{1},...,l_{i-1}%
,l_{i}+1,l_{i}+1,l_{i+1},...,l_{n})$, for any $1\leq i\leq n$, the
\textit{upper expansion} $l^{+}$ is the expansion in position $n$, while the
\textit{lower expansion} $l_{+}$ is the expansion in position $\max(1,$
$n-suf$ $l).$ Thus lower expansion is defined even for constant sequences, and
a sequence is constant if and only if its lower and upper expansion coincide.
Also $l\trianglelefteq h$ implies $l_{+}\trianglelefteq h_{+}$ in all cases,
just as it implies $l^{+}\trianglelefteq h^{+}$.

\bigskip

For path-length sequences with $n\geq2$ components and suffix length $k$, the
\textit{contraction} $\widehat{l}=(l_{1},...,l_{n-k},l_{n-k+1}-1,l_{n}%
,...,l_{n})$ is defined, it has $n-1$ components, and it is also a path-length
sequence, satisfying $\widehat{l}_{+}\trianglelefteq l\trianglelefteq
\widehat{l}^{+}$ .

\bigskip

\bigskip

\textbf{2 Uniqueness of the lattice meet}\bigskip

Stott Parker and Prasad Ram $\left[  \text{SPPR}\right]  $ state the fact that
the imbalance order on the set of path-length sequences with a given number of
components is a lattice. We show that this is can be verified as a consequence
of the Lemma below. It is of course enough to prove that the imbalance order
on path-length sequences $(l_{1},...,l_{n})$ with $n$ components is a meet semilattice.

\bigskip

Observe first that for any path-length sequences if $l\trianglelefteq h$ then
$last$ $l$ $\leq$ $last$ $h,$ and if $last$ $l$ $=last$ $h$ then $suf$ $l\leq
suf$ $h$ .

\bigskip

\textbf{Lemma}. \textit{For path-length sequences} $l\trianglelefteq h$, if
$last$ $l$ $<$ $last$ $h$ \textit{then} $l^{+}\trianglelefteq h_{+}$ . \ 

\bigskip

\textbf{Proof.} Let $l=$ $(l_{1},...,l_{n})$ and $h=$ $(h_{1},...,h_{n}).$ If
the suffix length of $h$ is $2k$, then for $i=n-2k$ we have $h_{1}\leq...\leq
h_{i}<h_{i+1}=...=h_{n}>l_{n}\geq...\geq l_{i+1}$ .

Recall that $l^{+}=(l_{1},...,l_{n-1},l_{n}+1,l_{n}+1)$ and $h_{+}%
=(h_{1},...,h_{i-1},h_{i}+1,h_{i}+1,h_{i+1},...,h_{n}).$

We have%

\[
2^{-h_{i}}\leq\Sigma\exp(h_{i+1},...,h_{n})\leq\frac{1}{2}\Sigma\exp
(l_{i+1},...,l_{n})\leq\Sigma\exp(l_{i+1},...,l_{n-k})
\]

It can be deduced that $\Sigma\exp(h_{1},...,h_{i-1},h_{i}+1)>\Sigma\exp
(h_{1},...,h_{i-1})\geq$ $\Sigma\exp(l_{1},..,l_{i-1},l_{i})$.

Also, for every $2\leq j\leq2k-1$, we have $\Sigma\exp(h_{n-j+1}%
,...,h_{n})\leq\Sigma\exp(l_{i-j+2},..,l_{n})=\Sigma\exp(l_{i-j+2}%
,...,l_{n}+1,l_{n}+1)$.

It follows that $l^{+}\trianglelefteq h_{+}.$
\ \ \ \ \ \ \ \ \ \ \ \ \ \ \ \ \ \ \ $\square$

\bigskip

Then let us prove by induction on $n$ the following:

\bigskip

\textbf{Proposition 1 \ }\textit{Path-length sequences} \textit{with the same
number} \textit{of components that is at most} $n$ \textit{alway}s
\textit{have a greatest lower bound (meet) in the imbalance order, for which
}$last$ $(s\wedge t)=\min(last$ $s,last$ $t)$\textit{.} \textit{For any such
sequences} $s$ \textit{and} $t$ \textit{with at least two components}
\textit{we have }$s\wedge t=(\widehat{s}\wedge\widehat{t})^{+}$\textit{\ if}
$(\widehat{s}\wedge\widehat{t})^{+}\trianglelefteq s,t$, \textit{otherwise
}$s\wedge t=(\widehat{s}\wedge\widehat{t})_{+}$ \textit{.}

\ 

\textbf{Proof.} \ The statement is obvious for $n=1.$ The inductive step from
$n-1$ to $n$ is as follows.

\bigskip

\textbf{Case 1}: $(\widehat{s}\wedge\widehat{t})^{+}\trianglelefteq s,t.$

\bigskip

To show that $(\widehat{s}\wedge\widehat{t})^{+}$ is the greatest lower bound
of $s,t$, let $l\trianglelefteq s,t$. Then $l\trianglelefteq s,t$ and thus
$\widehat{l}\trianglelefteq\widehat{s}\wedge\widehat{t}.$ From this
$l\trianglelefteq(\widehat{l})^{+}\trianglelefteq(\widehat{s}\wedge\widehat
{t})^{+}.$

\bigskip

From $(\widehat{s}\wedge\widehat{t})^{+}$ $\trianglelefteq s,t$ we get
\[
\min(last\text{ }\widehat{s},last\text{ }\widehat{t})+1\leq\min(last\text{
}s,last\text{ }t)
\]

\bigskip If $last$ $s=last$ $t$ then $suf$ $s$ or $suf$ \ $t$ is $2$, because
if both were larger, then $last$ $\widehat{s}$ $=$ $last$ $s$ and $last$
$\widehat{t}=last$ $t$, contradicting the above inequality.

But then $last$ $\widehat{s}+1=last$ $s$ or $last$ $\widehat{t}+1=last$ $t$,
and by the inductive hypothesis $last$ $(\widehat{s}\wedge\widehat{t})=(last $
$s)-1=(last$ $t)-1$ and $last$ $(\widehat{s}\wedge\widehat{t})^{+}=last$
$s=last$ $t=\min(last$ $s,last$ $t)$

\bigskip

If $last$ $s<last$ $t$ then $suf$ $s$ $=2$ because otherwise $last$
$\widehat{s}$ $=last$ $s\leq last$ $\widehat{t}$ and $(last$ $s)+1=$ $(last$
$\widehat{s})+1=$ $last$ $(\widehat{s}\wedge\widehat{t})^{+}\leq last$ $s$,
which is impossible. Also $last$ $\widehat{s}$ $<last$ $\widehat{t}.$

And also $(last$ $\widehat{s})$ $+1=last$ $s$ implying $\ last$ $(\widehat
{s}\wedge\widehat{t})^{+}=last$ $(\widehat{s}\wedge\widehat{t})+1=\min(last$
$\widehat{s},last$ $\widehat{t})+1=(last\widehat{s})+1=last$ $s=\min(last$
$s,last$ $t).$

\bigskip

\textbf{Case 2}: $(\widehat{s}\wedge\widehat{t})^{+}\ntrianglelefteq s$ or
$(\widehat{s}\wedge\widehat{t})^{+}\ntrianglelefteq t.$

\bigskip

Certainly still $(\widehat{s}\wedge\widehat{t})_{+}\trianglelefteq s,t$ . Note
that in this case $(\widehat{s}\wedge\widehat{t})$ cannot be constant.

By the induction hypothesis $last$ $(\widehat{s}\wedge\widehat{t})_{+}=last$
$(\widehat{s}\wedge\widehat{t})=\min(last$ $\widehat{s},last$ $\widehat{t}).$

\bigskip

\textbf{Subcase 2.1:} \ $last$ $s\neq last$ $t$ $\ $

Without loss of generality, we may suppose that $last$ $s<last$ $t.$

We claim that $suf$ $s>2.$ For if $suf$ $s=2$ then $s=(\widehat{s})^{+}$ and
$(\widehat{s}\wedge\widehat{t})^{+}\trianglelefteq s$ but $(\widehat{s}%
\wedge\widehat{t})^{+}\ntrianglelefteq t$. Also $last$ $\widehat{s}<last$
$\widehat{t}$, and thus $last$ $(\widehat{s}\wedge\widehat{t})=last$
$\widehat{s}$ is less then $last$ $\widehat{t}$, \ 

implying $(\widehat{s}\wedge\widehat{t})^{+}\trianglelefteq(\widehat{t}%
)_{+}\trianglelefteq t$ by the Lemma, a contradiction proving that $suf$ $s>2
$ .

Clearly then $last$ $(\widehat{s}\wedge\widehat{t})=\min(last$ $\widehat
{s},last$ $\widehat{t})=last$ $\widehat{s}=last$ $s.$

\bigskip

\textbf{Subcase 2.2:} \ $last$ $s=last$ $t$

Now the suffixes of both $s$ and $t$ cannot be $2$, because in that case
$(\widehat{s}\wedge\widehat{t})^{+}\trianglelefteq s,t$

If one of the suffix lengths, say $suf$ $s$ were $2,$ then $last$
$(\widehat{s}\wedge\widehat{t})=\min(last$ $\widehat{s},last$ $\widehat
{t})=last$ $\widehat{s}<last$ $\widehat{t}$. Then $(\widehat{s}\wedge
\widehat{t})^{+}\trianglelefteq(\widehat{s})^{+}=s$ and by the Lemma
$(\widehat{s}\wedge\widehat{t})^{+}\trianglelefteq(\widehat{t})_{+}=t,$

a contradiction.

Thus both suffix lengths are greater than $2,$ we have $(\widehat{s})_{+}=s$
and $(\widehat{t})_{+}=t$ and $last$ $(\widehat{s}\wedge\widehat{t}%
)=\min(last$ $\widehat{s},last$ $\widehat{t})=last$ $\widehat{s}=last$
$s=last$ $\widehat{t}=last$ $t$ .

\bigskip

In conclusion, in both subcases, if $last$ $s\leq last$ $t$ then $last$
$(\widehat{s}\wedge\widehat{t})=\min(last$ $\widehat{s},last$ $\widehat
{t})=last$ $\widehat{s}=last$ $s.$

\bigskip

We now return to the general conditions of Case 2. Again, without loss of
generality, we may suppose that $last$ $s\leq last$ $t.$

\bigskip

Let $l\trianglelefteq s,t.$ Obviously $\widehat{l}\trianglelefteq\widehat
{s}\wedge\widehat{t}$ and $last\widehat{l}\leq last$ $(\widehat{s}%
\wedge\widehat{t}).$

If $last\widehat{l}<last$ $(\widehat{s}\wedge\widehat{t})$ then by the Lemma
$l\trianglelefteq(\widehat{l})^{+}\trianglelefteq(\widehat{s}\wedge\widehat
{t})_{+}$ .

If $last\widehat{l}=last$ $(\widehat{s}\wedge\widehat{t})$ then $last$ $l\geq
last\widehat{l}=last$ $(\widehat{s}\wedge\widehat{t})=last\widehat{s}=last$
$s$ and, since $l\trianglelefteq s$ implies $last$ $l\leq last$ $s$, we have
$last$ $l=last$ $s=last\widehat{l}.$ This means that the suffix length of $l$
is also greater than $2$ and $l=(\widehat{l})_{+}\trianglelefteq(\widehat
{s}\wedge\widehat{t})_{+}$ \ \ \ \ \ \ \ \ $\square$

\bigskip

\bigskip

\textbf{3 Join-irreducible path-lengh sequences}

\bigskip

Stott Parker and Prasad Ram have shown $\left[  \text{SPPR}\right]  $ that for
path-length sequences the imbalance order relation $\trianglelefteq$ is the
transitive-reflexive closure the \textit{minimal balancing relation }which can
be defined as follows. For any path length sequence $l=$ $(l_{1},...,l_{n})$
call an index $1<j<n$ an \textit{excess} \textit{index} if $l_{j-1}%
<l_{j}=l_{j+1}$ and there is an index $i$ such that $l_{i}\leq l_{j}-2$. For
every excess index $j$ of $l$ consider the last index $i$ such that $l_{i}\leq
l_{j}-2,$ and let $bal\left[  l,j\right]  $ be the path-length sequence
obtained from $l$ by replacing the last occurrence of $l_{i}$ with two
consecutive entries equal to $l_{i}+1$ and replacing the first two
occcurrences of $l_{j}$ by a single entry equal to $l_{j}-1.$ The minimal
imbalance relation, on the set of path-length sequences $l$ with $n$
components, is
\[
\left\{  (\text{ }bal\left[  l,j\right]  ,\text{ }l\text{ }):\text{ }j\text{
is an excess index of }l\right\}
\]
The minimal imbalance relation contains (generally properly) the covering
relation of the partial order $\trianglelefteq$ .

\bigskip

Recall that an element of a finite lattice is\textit{ join-irreducible} if it
covers a unique element of the lattice (its \textit{unique} \textit{lower
cover}). The following is not difficult to verify:

\textbf{\bigskip}

\textbf{Proposition 2} \ \textit{Let} $j$ \textit{be the first excess index of
a path-length sequence} $l$. \textit{Then} $l$ \textit{is join-irreducible if
and only if for all excess indices} $k$ \textit{of} $l$ \textit{we have}
$bal\left[  l,j\right]  \trianglerighteq bal\left[  l,k\right]  .$

\bigskip

We call a sequence of integers \textit{near-constant} if it contains at most
two distinct values, and the difference between these is $1.$

\bigskip

\textbf{Proposition 3}\textit{ \ A path-length sequence } $l=$ $(l_{1}%
,...,l_{n})$ \textit{is join-irreducible if and only if it is the
concatenation of }$3$ \textit{sequences, }$l=uvw$ \textit{such that}

(i) \textit{both} $u$ and $w$ \textit{are near-constant, and} $v$ \textit{is
strictly increasing,}

(ii) $uv$ \textit{is not the empty sequence.}

(iii) \textit{if} $w$ \textit{is not constant but its first two components are
equal, than the value of these components is at least} $\left(  last\text{
}uv\right)  +2$.

\bigskip

\textbf{Proof} \ The conditions are easily seen to \ be sufficient for irreducibility.

Suppose on the other hand that $l$ is irrdeducible. It is useful to keep in
mind that there is a topological tree whose path-length sequence is $l$.
Obviously $l$ is not near-constant. Let $u$ be the longest near-constant
prefix of $l$, and $\ z$ its complement, $l=uz.$ Let $w$ be the longest
near-constant post-fix of $z$ and $v$ its complement, $z=vw.$ Clearly $w$
cannot be empty and, because $l=uvw$ is a path-length sequence, $w$ ends with
an even number of equal components, and the index of the first one of these is
an excess index $k.$ Also $uv$ cannot be empty, as required by condition (ii).

If $v$ had any repeated components, the index $j$ of the first of the first
two repeated components would be the first excess index of $l$. Then, by
Proposition 2, $l\vartriangleright bal\left[  l,j\right]  \vartriangleright
bal\left[  l,k\right]  $ would have to hold. In the tree corresponding to a
given path-length sequence, consider the number of nodes that are on some
root-to-leaf path of length at most $l_{j}$ This parameter is monotone, it can
never decrease as we go down in the imbalance lattice, but for $bal\left[
l,k\right]  $ it is the same as for $l$ itself. However, for $bal\left[
l,j\right]  $ the parameter is actually higher than for $l.$ This shows that
$v$ is strictly increasing, thus condition (i) also holds.

Finally, if condition (iii) failed, the first indices of the two constant runs
of components in $w,$ say $i$ and $k$ in that order, would both be excess
indices, and we would have $l_{i}=\left(  last\text{ }uv\right)  +1$. For the
first excess index $j$ we would have to have $j\leq i<k$ and thus
$l\vartriangleright bal\left[  l,j\right]  \vartriangleright bal\left[
l,k\right]  $ Consider now the sum of components of a path-length sequence.
This parameter, of integer value, is strictly monotone in the imbalance
lattice, it always decreases as we go down $\left[  \text{SPPR}\right]  $. But
from $l$ to $bal\left[  l,k\right]  $ it decreases only by $1$. This shows
that (iii) must also
hold.\ \ \ \ \ \ \ \ \ \ \ \ \ \ \ \ \ \ \ \ \ \ \ \ \ \ \ \ \ \ \ \ \ \ \ \ \ \ \ \ \ \ \ \ \ \ \ \ \ \ \ \ \ \ \ \ \ \ \ \ \ \ \ \ \ \ \ \ \ \ \ \ \ \ \ \ \ \ \ \ \ \ \ \ \ \ \ \ \ $\square
$

\bigskip

\bigskip

\textbf{Acknowledgements.}

This work has been co-funded by Marie Curie Actions and supported by the
National Development Agency (NDA) of Hungary and the Hungarian Scientific
Research Fund (OTKA, contract number 84593), within a project hosted by the
University of Miskolc, Department of Analysis.

\bigskip

\includegraphics[height=15mm, width=20mm]{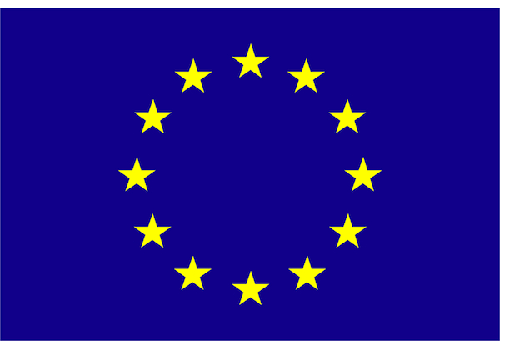}
\includegraphics[height=15mm, width=20mm]{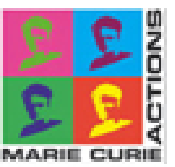}
\includegraphics[height=15mm, width=20mm]{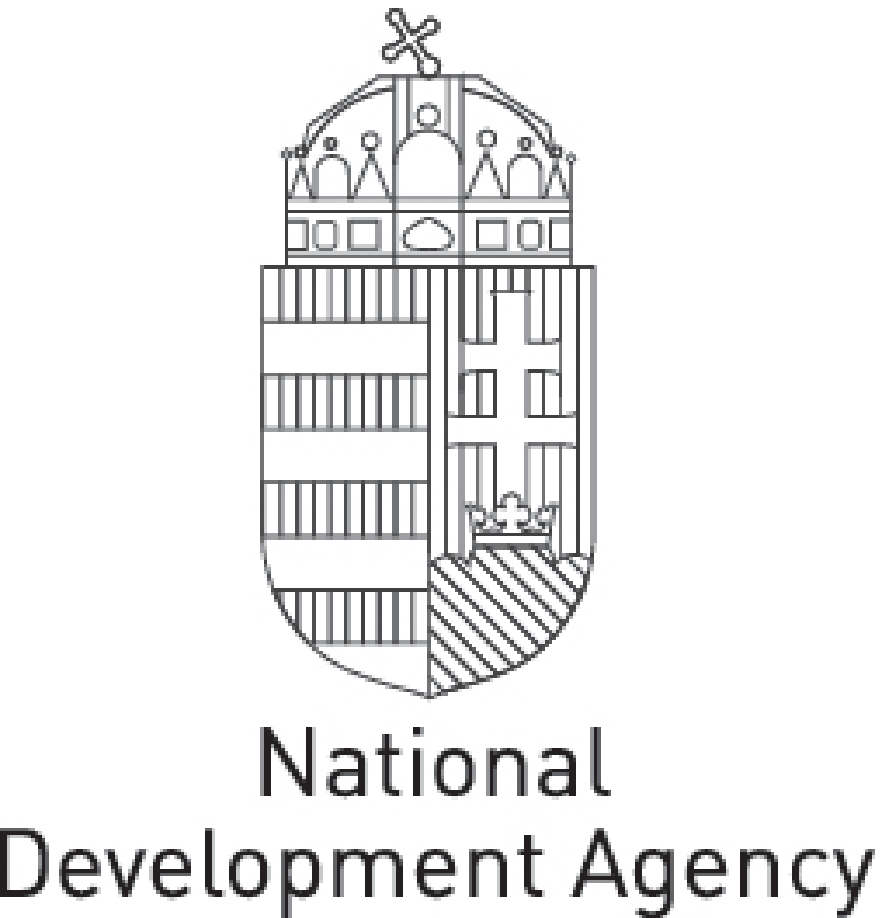} \includegraphics[height=15mm, width=20mm]{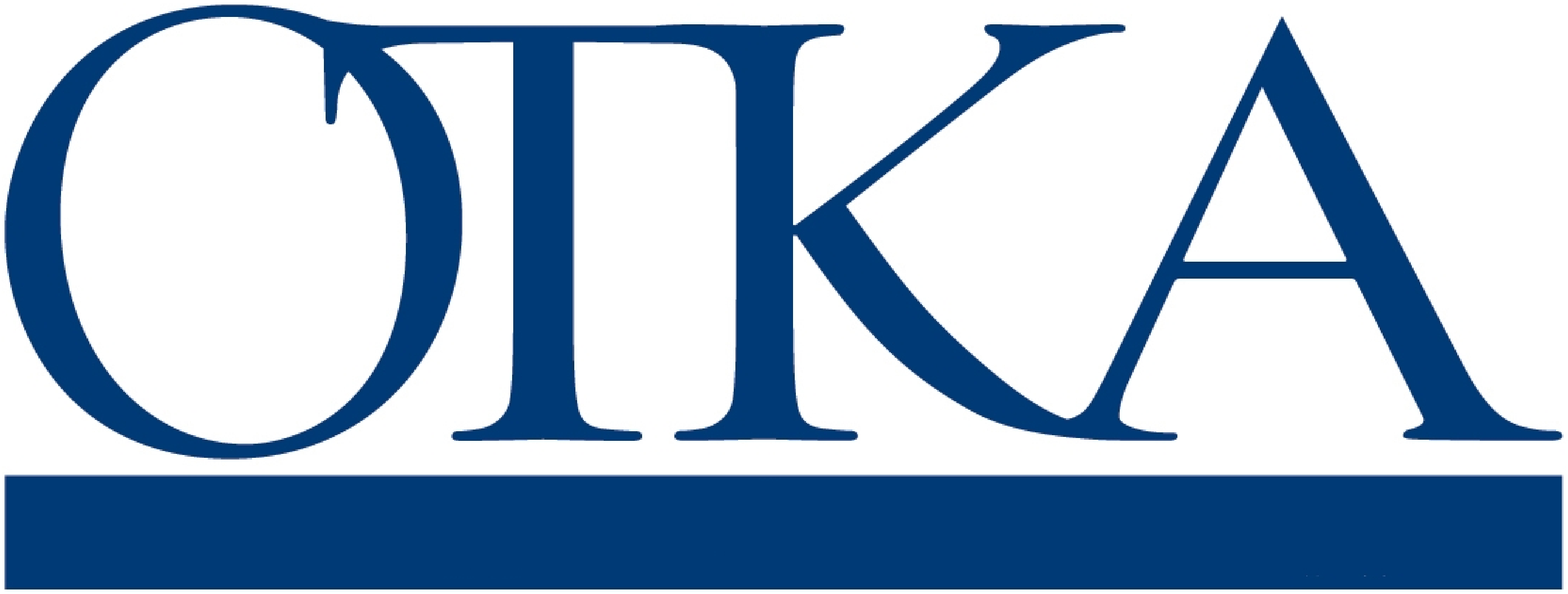}

\bigskip

\textbf{References}

\bigskip

$\left[  \text{K}\right]  $ \ L.G. Kraft, \textit{A Device for Quantizing,
Grouping, and Coding Amplitude Modulated Pulses}, Q.S. Thesis, MIT 1949

\bigskip

$\left[  \text{SPPR}\right]  $ \ D. Stott Parker, Prasad Ram, The Construction
of Huffman Codes is a Submodular ("Convex") Optimization Problem Over a
Lattice of Binary Trees. \textit{SIAM J. Comput.} 28(5) 1875-1905 (1999)

\end{document}